\documentclass{ifacconf}

\usepackage{graphicx}      
\usepackage{natbib}        
\usepackage{amssymb, amsmath}
\begin{document}
\begin{frontmatter}

\title{Energy Estimates for Low Regularity Bilinear Schr\"{o}dinger Equations 
\thanksref{footnoteinfo}} 

\thanks[footnoteinfo]{This work has been partially supported by INRIA Nancy-Grand 
Est. Second and third authors were partially supported by
French Agence National de la Recherche ANR ``GCM'' program
``BLANC-CSD'', contract number NT09-504590. The
third author was partially supported by European Research
Council ERC StG 2009 ``GeCoMethods'', contract number
239748.}

\author[First]{Nabile Boussa\"{i}d} 
\author[Second]{Marco Caponigro} 
\author[Third]{Thomas Chambrion}

\address[First]{Laboratoire de math\'ematiques, Universit\'e de Franche--Comt\'e,
25030 Besan\c{c}on, France (e-mail:  Nabile.Boussaid@univ-fcomte.fr)}
\address[Second]{D\'epartement  Ing\'enierie Math\'ematiques, Conservatoire National des Arts et M\'eti\'ers, 75003 Paris, France (e-mail:
marco.caponigro@cnam.fr)}
\address[Third]{Universit\'e de Lorraine, Institut \'Elie Cartan de Lorraine, 
Vand{\oe}uvre-l\`es-Nancy, F-54506, France,\\  CNRS, IECL, 
Vand{\oe}uvre-l\`es-Nancy, F-54506, France,\\ Inria, CORIDA, Villers-l\`es-Nancy, F-54600, France (e-
mail:Thomas.Chambrion@univ-lorraine.fr)}

\begin{abstract}                
This paper presents an energy estimate in terms of the total variation of the 
control for bilinear infinite dimensional quantum systems with unbounded 
potentials. These estimates allow a rigorous construction of  
propagators associated with  controls of bounded variation. Moreover, upper bounds of the error 
made when replacing the infinite dimensional system by its finite dimensional Galerkin
approximations is presented. 
\end{abstract}


\end{frontmatter}

\section{Introduction}

\subsection{Physical context}
The state of a quantum system evolving in a Riemannian manifold $\Omega$ is
described by 
its \emph{wave function}, a point $\psi$ in $L^2(\Omega, \mathbf{C})$. When the
system is 
submitted to an electric field (e.g., a laser), the time evolution of the wave
function is 
given, under the dipolar approximation and neglecting decoherence,  by the
Schr\"{o}dinger 
bilinear 
equation:
\begin{equation}\label{EQ_bilinear}
\mathrm{i} \frac{\partial \psi}{\partial t}=(-\Delta  +V(x)) \psi(x,t) +u(t) 
W(x) 
\psi(x,t)
\end{equation}
where $\Delta$ is the Laplace-Beltrami operator on $\Omega$,  $V$ and $W$ are
real 
potential accounting for the properties of the free system and the control
field 
respectively, while 
the real function of the time $u$ accounts for the intensity of the laser. 

In view of applications (for instance in NMR), it is important to know whether and how
it is 
possible to chose a suitable control $u:[0,T]\to \mathbf{R}$ in order to steer 
(\ref{EQ_bilinear}) from a 
given initial state to a given target. This question has raised considerable
interest in 
the community in the last decade. After the negative results of \cite{bms} and 
\cite{Turinici} excluding 
exact controllability on the natural domain of the operator $-\Delta +V$ when
$W$ is bounded, the first, and 
at this day the only one, description of the attainable set for an example of
bilinear 
quantum system 
was obtained by (\cite{beauchard,beauchard-coron}).  Further investigations of
the 
approximate controllability of (\ref{EQ_bilinear}) were conducted using
Lyapunov 
techniques  
(\cite{nersesyan, Nersy, beauchard-nersesyan,Mirrahimi, MR2168664, mirra-solo}) 
and 
geometric techniques (\cite{Schrod,Schrod2}).

In  most of the references cited above, the potentials $V$ and $W$ in
(\ref{EQ_bilinear}) 
are bounded. The very general (and irregular) systems considered by
\cite{Schrod2} allow 
to define the solutions of (\ref{EQ_bilinear}) for piecewise constant controls
only. The 
aim of this paper is to present a coherent framework  to deal with unbounded
potentials in 
(\ref{EQ_bilinear}).  This includes a rigorous definition of the solution of 
(\ref{EQ_bilinear}) for control that are not necessarily piecewise constant and
the 
extension of some quantitative energy estimates.

\subsection{Abstract framework and notations}
We reformulate the control problem in more abstract framework, in such a way
that we can use some of the  powerful tools of functional analysis. In a
separable Hilbert space $H$, we consider a pair $(A,B)$ of (possibly unbounded)
linear operators  that satisfy Hypothesis \ref{HYP_1}
\begin{hypo}\label{HYP_1}
$(A,B)$ is a pair of linear operators such that
\begin{enumerate}
\item $A$ is skew-adjoint on its domain $D(A)$; \label{ASS_A_skew_adjoint}
 \item $\mathrm{i}A$ is bounded from below; \label{ASS_iA_bounded_from_below}
\item $B$ is skew-symmetric;  \label{ASS_B_skew_sym}
\item there exists $a,b\geq 0$ such that $\|B\psi\| \leq a\|A \psi\| +b
\|\psi\|$ for any $\psi$ 
in $D(A)$.\label{ASS_B_relativ_bounded}
\end{enumerate}
\end{hypo}
Following \cite{Kato1953}, Hypothesis \ref{HYP_1} is the minimal framework
for our developments.
For many examples encountered in the physics literature, $A$ has a discrete
spectrum
and we will consider the more restrictive Hypothesis \ref{HYP_2}.
\begin{hypo}\label{HYP_2}
$(A,B,(\phi_j)_{j\in \mathbf{N}},\alpha)$ is a quadruple such that
\begin{enumerate}
\item $(A,B)$ satisfies Hypothesis \ref{HYP_1}; 
\item $(\phi_j)_{j\in \mathbf{N}}$ is a Hilbert basis of
$H$;\label{ASS_base_Hilbert}
\item $0\leq \alpha \leq1$;
\item $A$ has discrete spectrum $(-\mathrm{i}\lambda_j)_{j \in \mathbf{N}}$
with 
$\lambda_j\to +\infty$ as $j \to \infty$;\label{ASS_discrete_spectrum}
\item for any $j$ in $\mathbf{N}$, $A\phi_j=-\mathrm{i} \lambda_j \phi_j$; 
\label{ASS_vect_propres}
\item there exists $d\geq 0$ such that $\|B\psi\| \leq d\||A|^\alpha \psi\|$ for
any  $\psi$ 
in $D(|A|^\alpha)$.\label{ass:bdd}
\end{enumerate}
\end{hypo}

Thanks to the Kato-Rellich theorem (see \cite{kato}), Hypotheses 
\ref{HYP_1}.\ref{ASS_A_skew_adjoint}, \ref{HYP_1}.\ref{ASS_B_skew_sym} and
\ref{HYP_1}.\ref{ASS_B_relativ_bounded} imply that, for any $u$ in
$(-1/a,1/a)$, $A+uB$ is skew-adjoint with domain $D(A)$ and generates a unitary
propagator $t\mapsto e^{t(A+uB)}$. In  
particular, this allows to define by concatenation the propagator $\Upsilon^u:t
\mapsto 
\Upsilon^u_t$ for the control system
\begin{equation}\label{EQ_main}
\frac{\mathrm{d}\psi}{\mathrm{d}t}=(A+u(t)B)\psi
\end{equation}
for $u$ piecewise constant $u$ taking value in $(-1/a,1/a)$.

Recall that a function $u:[0,T]\to \mathbf{R}$ has \emph{bounded variation}
(or is BV) if
there exists a constant $C$ such that, for any partition $0=a_0<a_1<\ldots
<a_n=T$ of $[0,T]$, $\sum_{k=1}^n |u(a_k)-u(a_{k-1)})|<C$. The 
smallest $C$  satisfying this property for any partition of $[0,T]$ is the
\emph{total variation} of  $u$, denoted  $TV_{[0,T]}(u)$. 

We define the set $\mathcal{U}$ of the functions $u:\mathbf{R}\to \mathbf{R}$
with  
bounded variation such 
that $u(t)=0$ for $t\leq 0$. In $\mathcal U$, the sequence $(u_n)_{n\in
\mathbf{N}}$ 
converges to $u$ if 
$\sup_n  TV_{\mathbf{R}}(u_n)\leq TV_{\mathbf{R}}(u)$ and $u_n(t)$ tends to
$u(t)$ as $n$ 
goes to infinity for almost any $t$ in 
$\mathbf{R}$. 

\subsection{Contribution of this paper}

This paper presents a rigorous yet elementary construction 
of the solutions of (\ref{EQ_main}) associated with controls of bounded variation,
inspired from \cite{Kato1953}. Among other byproducts of our energy estimates,
we give a lower bound for the number of switches needed to steer (\ref{EQ_main})
from a 
given source to a 
given target using controls with value in $\{0,1\}$ and we give an upper bound
of the error made when one 
replaces the original infinite dimensional system (\ref{EQ_main}) by one of its
finite dimensional Galerkin approximation.
Such estimates are instrumental in practice, both for theoretical analysis,
design of control laws and numerical simulations. 

The strength of our results is the relative generality of our assumptions. In
this sense, this paper may be seen as an extension of the results
of \cite{ACCFEPS} to systems that are not \emph{weakly-coupled} (according to \cite[Definition~1]{weakly-coupled}).

\subsection{Content of the paper}

The first part of the paper (Section \ref{SEC_construction_propagators})
 is concerned with the construction of the 
solutions of (\ref{EQ_main}) for controls with bounded variation.
The key point of this construction is an energy estimate in terms of the total
variation
of the control (see Proposition \ref{PRO_Borne_A}). 
The second part of the paper (Section \ref{SEC_GGA}) presents some consequences
of this energy estimate
in terms of approximation of the original infinite dimensional system by its
finite dimensional dynamics.
Finally, we apply our results to various types of quantum oscillators 
encountered in the physics literature (Section \ref{SEC_examples}).

\section{Construction of the propagators}\label{SEC_construction_propagators}

To begin with, we consider the simple case where $\|B\psi\|\leq a
\|A\psi\|$ for 
any $\psi$ in $D(A)$. The general case of operators $B$ relatively bounded
with respect 
to $A$ satisfying Hypothesis \ref{HYP_1}.\ref{ASS_B_relativ_bounded} will be 
treated in Subsection 
\ref{SEC_constr_B_general}. 

\subsection{Estimates on the $A$ norm}

For any $\psi$ in $D(A)$, for any $u$ in $\mathbf{R}$ such that
$a|u|<1$, 
\begin{eqnarray}
\|B\psi\| & \leq &a \|A \psi \| \\
&\leq & a (\|(A+uB)\psi \| + |u| \|B \psi \|)
\end{eqnarray}  
Hence,
\begin{eqnarray}
(1-a |u|) \|B\psi\|&\leq& a \|(A+u B)\psi \|\\
\|B \psi \| &\leq& \frac{ a}{1-|u| a}\|(A +uB)\psi \|
\end{eqnarray}  

For any $u_1,u_2$ in $(-1/a,1/a)$, $t$ in $\mathbf{R}$ and $\psi$
in $D(A)$, $\psi$ is
in $D(A+u_2B)$ by Hypothesis \ref{HYP_1}.\ref{ASS_B_relativ_bounded}. Hence,
$e^{t(A+u_2B)} \psi$ belongs to 
$D(A+u_2B)=D(A)=D(A+u_1B)$. Moreover,
\begin{eqnarray*}
\lefteqn{\|(A+u_1B) e^{t(A+u_2B)}\psi \|} \nonumber\\
&\leq & \| (A+u_2B) e^{t(A+u_2B)}\psi \| + \|(u_1-u_2)B e^{t(A+u_2B)}\psi \| \\
& \leq & \| e^{t(A+u_2B)}(A+u_2B) \psi  \|  \nonumber \\
& & \quad + |u_1-u_2|  \frac{a}{1-|u_2| a} \|(A
+u_2B)e^{t(A+u_2B)}\psi \|
\end{eqnarray*} 
and hence
\begin{align*}
\|(A+u_1B) e^{t(A+u_2B)}\psi \|&\\ \leq \left (1+\frac{|u_1-u_2|
\|B\|_{A}}{1-|u_1| a} \right )  &\|(A+u_2B)\psi \|
\end{align*}
and
\begin{eqnarray*}
\|(A+u_1B) \psi \|\leq \left (1+\frac{|u_1-u_2|
a}{1-|u_1| a} \right )  \|(A+u_2B)\psi \|
\end{eqnarray*}

Let $u^\ast>0$ be given such that $  |u^\ast| <1/a$. For any $t\geq 0$,
for any $u_1,u_2$ in $(-u^\ast,u^\ast)$ one has, with
$\Gamma=\frac{a}{1- |u^\ast| a}$, 
\begin{eqnarray*}
\lefteqn{\|(A+u_1B) e^{t(A+u_2B)}\psi \|}\\
& \quad \quad \quad  \leq& \exp \left ( \Gamma |u_2-u_1|  \right) \|(A+u_2B)
\psi \|.
 \end{eqnarray*}

Consider now a piecewise constant control $u:[0,T]\to (-1/a,1/a)$
taking value $u_j$ for
time $t_j$, $t_j\geq 0$ $1\leq j \leq p$, $p\in \mathbf{N}$. 
We get by concatenation, for any $\psi$ in $D(A)$,  
\begin{eqnarray*}
\lefteqn{\|A \Upsilon^u_{T,0}\psi \|}\\
 & \leq & \exp(\Gamma |u_p|) \times \\
 &&\times\|(A+u_p B)e^{t_p(A+u_pB)}
e^{t_{p-1}(A+u_{p-1}B)}\cdots e^{t_1(A+u_1B)}\psi\|\\
&\leq & \exp(\Gamma |u_p|) \left \lbrack \prod_{j=1}^p \exp \left
( \Gamma
|u_j-u_{j+1}|  \right) \right \rbrack \exp(\Gamma |u_1|) \|A \psi \|\\
& \leq & \exp (2 \Gamma TV_{[0,T]}(u))  \|A \psi \|.
\end{eqnarray*}

We obtain, similarly to \cite{Kato1953}, the following result.
\begin{prop}\label{PRO_Borne_A}
For any $\delta\in(0,1)$, let $(A,B)$ satisfy Hypothesis \ref{HYP_1}. Then,
for
any piecewise constant $u:[0,T]\to
(-(1-\delta)/a,(1-\delta)/a)$, for any $\psi$ in $D(A)$,
$\displaystyle{
\|A \Upsilon^u_{T,0}\psi\|\leq e^{\frac{2a}{\delta} TV_{[0,T]}(u)} \|A \psi
\|}$.
\end{prop}

\subsection{Definition of propagators for BV controls}

For any $\delta\in(0,1)$ and $a > 0$, let $\mathcal{U}_{\delta,a}$ be the subset of
$u\in \mathcal{U}$ such that
$u:\mathbf{R} \to (-(1-\delta)/a,(1-\delta)/a)$.

Let $u$ in $\mathcal{U}_{\delta,a}$.
There exists a sequence $u_n$ in $\mathcal{U}_{\delta,a}$ of piecewise constant
functions such that 
(i) $(u_n)_n$ tends to $u$ pointwise and
(ii) for any $n$ in $\mathbf{N}$, $TV_{[0,T]}(u_n)\leq TV_{[0,T]}(u)$. 
These conditions implies that $\sup_n \|u_n\|_{L^\infty}<+\infty$. 

\begin{prop} Let $(A,B)$ satisfy Hypothesis~\ref{HYP_1} with $b=0$ and let $(u_n)_n$ be defined
as above.  
For any $t$ in $[0,T]$, for any $\psi$ in $D(A)$,
$(\Upsilon^{u_n}_{(t,0)}\psi)_{n \in \mathbf{N}}$ is a Cauchy sequence (for the norm of $H$).
\end{prop}
\begin{pf}
 By Duhamel's identity, for any $\psi$ in $D(A)$,
\begin{eqnarray*}
\Upsilon^{u_n}_{(t,0)}\psi-\Upsilon^{u_m}_{(t,0)}\psi= \int_0^t
\Upsilon^{u_n}_{(s,t)}(u_n(s)-u_m(s))B \Upsilon^{u_m}_{(s,0)}\psi\mathrm{d}s 
\end{eqnarray*}
For any $s$ in $(0,t)$, by Proposition~\ref{PRO_Borne_A}, 
$$\sup_{0\leq s\leq t \leq T}\sup_{n,m}\|\Upsilon^{u_n}_{(s,t)}B
\Upsilon^{u_m}_{(s,0)}\psi\|<+\infty. $$ Moreover, 
$u_n(s)-u_m(s)$ tends to zero as $n,m$  tend to infinity ($(u_l(s))_l$ is a
Cauchy sequence). The result follows from Lebesgue's dominated convergence
theorem.
\end{pf}

We define $\Upsilon^u_{(t,0)}\psi=\lim_n \Upsilon^{u_n}_{(t,0)}\psi$ for any
$\psi$ in $D(A)$. It is clear from the definition that
the construction is independent on the choice the sequence $(u_n)_n$
converging to $u$. 
Since $D(A)$ is dense in $H$ and $\Upsilon^u_{(t,0)}$ is bounded (in $H$ norm)
by 1 on $D(A)$, $\Upsilon^u_{(t,0)}$ admits an extension to $H$ that we still
denote with
 $\Upsilon^u_{(t,0)}$.

\subsection{General case of $A$-bounded operators} \label{SEC_constr_B_general}

Next proposition states that replacing $A$ by $A_\lambda:=A+\mathrm{i}\lambda \mathrm{Id}$ 
induces just a 
global phase shift at the level of the propagators.
\begin{prop}
 For any $\delta\in(0,1)$, any $u$ in
$\mathcal{U}_{\delta,a}$, for any $(A,B)$
satisfying 
Hypothesis~\ref{HYP_1} with $b=0$ in 
Hypothesis~\ref{HYP_1}.\ref{ASS_B_relativ_bounded},  for 
any $\lambda$ in $\mathbf{R}$, 
denote with $\Upsilon^{u}_{t,0}$ and $\Upsilon^{u,\lambda}_{t,0}$ the
propagators 
associated with $x'=(A+uB)x$ and $x'=(A_\lambda+uB)x$ respectively. Then
$\Upsilon^{u,
\lambda}_{t,0}=e^{\mathrm{i}\lambda t}\Upsilon^{u}_{t,0}$. 
\end{prop}
\begin{pf}
The result is obvious with piecewise constant controls. The result follows by
taking the limit for a sequence of piecewise constant controls $(u_n)_n$
tending to $u$ for the $BV$ topology. 
\end{pf}

We now come back to the definition of propagators of (\ref{EQ_main}) in the
general case $\|B\psi\| \leq a\|A \psi\|+b \|\psi\|$. As $A$ is bounded from below
(Hypothesis~\ref{HYP_1}.\ref{ASS_iA_bounded_from_below}), for every $\eta>0$, there exists $\lambda$ large
enough  such that $\|B\psi\| \leq
(a+\eta)\|A_\lambda \psi\|$ and we apply the 
above procedure to $(A_\lambda,B)$ to define the propagator
$\Upsilon^{u,\lambda}_{t,0}$ and, finally, the propagator
$\Upsilon^{u}_{t,0}:=e^{-\mathrm{i}\lambda
t}\Upsilon^{u,\lambda}_{t,0}$. Notice 
that this construction is independent on $\lambda$, provided that $\lambda$ is large
enough. Notice also, and this is instrumental in our study, that for any $u$ with bounded variation such that $\sup |u|<1/a$, for any $\lambda$ large enough,  $\|A \Upsilon^{u}_{t,0} \psi_0\| =\|A \Upsilon^{u,\lambda}_{t,0} \psi_0\|$ for every $\psi_0$ in $D(|A|)$. 

Below we write $\Upsilon^{u,\lambda}_{t}$ for $\Upsilon^{u,\lambda}_{t,0}$
and $\|\psi \|_r$ for $\|(1+|A|)^r\psi \|$.

We sum up the result of Section \ref{SEC_construction_propagators} in the following Proposition.
%
\begin{prop}\label{PRO_continuity_norme_Ar}
Let $\delta\in (0,1)$ and $(A,B)$ satisfy Hypothesis \ref{HYP_1}.
For any $u$ in
$\mathcal{U}_{\delta,a}$, for any $t\geq 0$, the 
propagator$\Upsilon^u_{t,0}:\psi \mapsto \Upsilon^u_{t,0} \psi  $ is
continuous from $D(|
A|)$ to $D(|A|)$. 
Moreover, for every $\eta>0$, there exists $\lambda \in \mathbf{R}$ such that, for any $\psi$ in $D(A)$, for any $u$ in $\mathcal{U}_{\delta,a}$, for any $t\geq 0$,
$$\|(A+\mathrm{i}\lambda) \Upsilon^u_{t,0}\psi\|\leq e^{\frac {2a+\eta}{\delta} TV_{[0,t]}(u)}
\|(A+\mathrm{i}\lambda) \psi \|.$$
\end{prop}

\section{Good Galerkin Approximations}\label{SEC_GGA}

For applications (design of control laws or numerical simulations), it is
common 
to replace the original infinite dimensional system (\ref{EQ_main}) by a suitable finite 
dimensional approximation. It is often possible to bound the error due to this
approximation.
Under Hypothesis \ref{HYP_2}, we derive in this section an explicit upper bound 
of this error that  depends only on the $L^1$ norm and the total variation of
the control. 
 The results presented here extend the results of \cite{ACCFEPS}. 

\subsection{Notion of Good Galerkin Approximations}
Let $\mathbf{\Phi}=(\phi_j)_{j \in \mathbf{N}}$ be a Hilbert basis of $H$.
For any $N$ in $\mathbf{N}$, we define the orthogonal projection
$$
 \pi_N^{\mathbf{\Phi}}\psi \in H \mapsto \sum_{j\leq N} \langle
\phi_j,\psi\rangle
\phi_j \in H.
$$
\begin{defn}\label{DEF_Galerkin_approx}
Let $(A,B,\Phi,1)$ satisfy Hypothesis \ref{HYP_2} and $N \in \mathbf{N}$.  The
\emph{Galerkin approximation}  of (\ref{EQ_main})
of order $N$ is the system in $H$ 
\begin{equation}\label{eq:sigma}
\dot x = (A^{(\mathbf{\Phi},N)} + u(t) B^{(\mathbf{\Phi},N)}) x 
\end{equation}
where {$A^{(\mathbf{\Phi},N)}= \pi_N^{\mathbf{\Phi}} A_{\upharpoonright_{\rm Im}
\pi_N^{\mathbf{\Phi}}}$ and
$B^{(\mathbf{\Phi},N)}= \pi_N^{\mathbf{\Phi}} B_{\upharpoonright_{\rm Im}
\pi_N^{\mathbf{\Phi}}}$} are the
\emph{compressions} of $A$ and $B$ (respectively).
\end{defn}

{We denote by $X^{u}_{(\mathbf{\Phi},N)}(t,s)$ the propagator of
(\ref{eq:sigma})
associated with a $L^{1}$ function $u$.}

\begin{rem}
The operators $A^{(\mathbf{\Phi},N)}$ and $B^{(\mathbf{\Phi},N)}$ are defined
on the \emph{infinite} dimensional space $H$. However, they have finite rank and
the dynamics of $(\Sigma_N)$ leaves invariant the $N$-dimensional space
$\mathcal{L}_{N} = \mathrm{span}_{1\leq j \leq N} \{\phi_j\}$. Thus,
$(\Sigma_N)$ can be seen as a finite dimensional bilinear system in
$\mathcal{L}_{N}.$
\end{rem}

The system $(A,B)$ admits a sequence of \emph{Good Galerkin Approximations} (GGA
in short), in time $T\in (0,+\infty]$, for a functional 
norm $N(\cdot)$ on a functional space $\mathbf{U}$ in a subspace  $D$ (with
norm $\|\cdot \|_D$) of $H$ if, 
for any $K,\varepsilon>0$, for any $\psi$ in $D$, there exists $N$ in
$\mathbf{N}$ such that,  for any $u$ 
in $\mathbf U$, $N(u)\leq K$ implies
$\|(X^{u}_{(\mathbf{\Phi},N)}(t,0)-\Upsilon^u_{t,0})\psi \|_D
<\varepsilon$ for any $t<T$.
 
\subsection{GGA for BV controls} 

\begin{prop}\label{PRO_troncature}
Let $(A,B,\Phi)$ satisfy
Hypotheses~\ref{HYP_1}, \ref{HYP_2}.\ref{ASS_base_Hilbert},
\ref{HYP_2}.\ref{ASS_discrete_spectrum} and \ref{HYP_2}.\ref{ASS_vect_propres}.
Then, for any $\delta\in(0,1)$, for any $r\in [0,1)$
for any  $n\in \mathbf{N}$, $N \in \mathbf{N}$,
$(\psi_j)_{1\leq j \leq n}$ in $D(|A|)^n$,
and for any function $u$ in $\mathcal{U}_{\delta,a}$,
\begin{equation}\label{eq:feps2}
\|(\mathrm{Id} - \pi^{\Phi}_{N})
\Upsilon^{u}_{t}(\psi_{j})\|_r\leq   \frac{e^{\frac 2\delta a
TV_{\mathbf{R}}(u)}\|\psi_j\|_{1}}{{\inf_{j>N}\lambda_{j}^{1-r}}}.
\end{equation}
for any $t \geq 0$ and $j=1,\ldots,n$.
\end{prop}

\begin{pf}
Fix $j \in \{1,\ldots,n\}$. For any $N > 1$, one has
\begin{eqnarray*}\label{eq:estimates}
\left\|(\mathrm{Id} - \pi^{\Phi}_{N})
\Upsilon^{u}_{t,0}(\psi_{j})\right\|_r^2&=&
\sum_{n = N+1}^{\infty}  \lambda_n^{2r}| \langle  \phi_{n},
\Upsilon^{u}_{t}(\psi_{j}) \rangle|^{2}\\
 &\leq&\inf_{j>N}\lambda_{j}^{2(r-1)}
\left\|\Upsilon^{u}_{t,0}(\psi_{j})\right\|_{1}^2.
\end{eqnarray*}
By Proposition~\ref{PRO_Borne_A}, for
any $t>0$,
$$\| \Upsilon^u_{t,0}\psi_j\|_1\leq e^{\frac 2\delta a
TV_{\mathbf{R}}(u)}
\|\psi_j\|_1.$$
\end{pf}

\begin{prop}[Good Galerkin Approximation]\label{PRO_Good_Galerkin_approximtion}
Let $\delta\in (0,1)$, $\alpha\in [0,1)$ and $(A,B,\Phi,\alpha)$ satisfy
Hypothesis~\ref{HYP_2}. Then
for any $\varepsilon > 0 $, $K\geq 0$, $n\in \mathbf{N}$, and
$(\psi_j)_{1\leq j \leq n}$ in $D(|A|)^n$
there exists $N \in \mathbf{N}$
such that
for any $L^1$ function $u$ in $\mathcal{U}_{\delta,a}$,
$$
\|u\|_{L^{1}} +TV_{\mathbf{R}}(u)< K \Rightarrow \| \Upsilon^{u}_{t}(\psi_{j}) -
X^{u}_{(\mathbf{\Phi},N)}(t,0)\pi_{N} \psi_{j}\| < \varepsilon,
$$
for any $t \geq 0$ and $j=1,\ldots,n$.
\end{prop}

\begin{pf}
Fix $j$ in $\{1,\ldots,n\}$ and consider
the map $t\mapsto
\pi_{N} \Upsilon^{u}_{t}(\psi_{j})$ that is absolutely continuous and satisfies,
for almost any $t \geq 0$,
\begin{eqnarray*}
\lefteqn{\frac{d}{dt} \pi_{N} \Upsilon^{u}_{t}(\psi_{j}) 
= (A^{(\mathbf{\Phi},N)} +  u(t)B^{(\mathbf{\Phi},N)}) \pi^{\Phi}_{N}
\Upsilon^{u}_{t}(\psi_{j})} \\
&&\quad \quad \quad \quad \quad \quad \quad \quad \quad \quad \quad \quad+ u(t)
\pi^{\Phi}_{N} B (\mathrm{Id} - \pi^{\Phi}_{N}) \Upsilon^{u}_{t}(\psi_{j}).
\end{eqnarray*}
Hence, by variation of constants, for any $t \geq 0$,
\begin{eqnarray}\label{EQ_preuve_good_Galerkin}
\lefteqn{\pi_{N} \Upsilon^{u}_{t}(\psi_{j})= X^{u}_{(\mathbf{\Phi},N)}(t,0) 
\pi^{\Phi}_{N}\psi_{j}} \nonumber \\  
&&\quad \quad +
 \int_{0}^{t} \!\!\!
X^{u}_{(\mathbf{\Phi},N)}(t,s) \pi^{\Phi}_{N} B(\mathrm{Id} - \pi_{N})
\Upsilon^{u}_{s}(\psi_{j})
u(\tau)  \mathrm{d}\tau.
\end{eqnarray}
By Proposition~\ref{PRO_troncature}, the norm of $t \mapsto B(\mathrm{Id} -
\pi_{N}) \Upsilon^{u}_{t}(\psi_{j})$ is
less than
$d e^{\frac2\delta a K}\inf_{j>N}\lambda_{j}^{\alpha-1}\|\psi_{j}
\|_1$. Since $X^{u}_{(\mathbf{\Phi},N)}(t,s) $ is unitary,
\begin{align*}
 \|\pi_{N} \Upsilon^{u}_{t}(\psi_{j}) -X^{u}_{(\mathbf{\Phi},N)}(t,0)&
\pi_{N}\psi_{j}\|\\ \leq
\|u\|_{L^1}&
d\inf_{j>N}\lambda_{j}^{\alpha-1}e^{c(A,B)K}\|\psi_{j}\|_1.
\end{align*}
Then
\begin{eqnarray*}
\lefteqn{\|\Upsilon^u_t(\psi_j) - X^u_{(N)}(t,0)\pi^{\Phi}_{N} \psi_{j}\|
 }\\
 &\quad \leq & \|(\mathrm{Id} -\! \pi_{N}) \Upsilon^{u}_{t}(\psi_{j})\| \! + \!
\|\pi^{\Phi}_{N} \Upsilon^{u}_{t}(\psi_{j}) \! - \! X^{u}_{(\mathbf{\Phi},N)}(t,0) 
\pi^{\Phi}_{N} \psi_{j}
\| \nonumber\\
& \quad \leq& \|u\|_{L^1}
(1+dK)d e^{\frac2\delta a K}\inf_{j>N}\lambda_{j}^{\alpha-1}\|\psi_{j}\|_1.
\label{eq:qwer}
\end{eqnarray*}
This completes the proof  since $\lambda_n$ tends to infinity as $n$
goes to infinity. \end{pf}

\section{Examples}\label{SEC_examples}

\subsection{Tri-diagonal systems}
\begin{defn}
A system $(A,B,\mathbf{\Phi})$ is \emph{tri-diagonal} if $(A,B)$ satisfies 
Hypotheses~\ref{HYP_1}.\ref{ASS_A_skew_adjoint}, 
\ref{HYP_1}.\ref{ASS_iA_bounded_from_below},
\ref{HYP_1}.\ref{ASS_B_skew_sym}, \ref{HYP_2}.\ref{ASS_base_Hilbert}, 
\ref{HYP_2}.\ref{ASS_discrete_spectrum} and \ref{HYP_2}.\ref{ASS_vect_propres}
and  if,
for any $j,k$ in $\mathbf{N}$,
$|j-k|>1$ implies $\langle \phi_j, B\phi_k \rangle =0$.
\end{defn}
In the following, we denote $b_{j,k}=\langle \phi_j, B\phi_k \rangle$.

\begin{prop}\label{PRO_diagonal_B_A_borne}
Let $(A,B,\mathbf{\Phi})$ be a tri-diagonal system and let $r$ be a positive
number. Assume
that the sequences 
$\left (\frac{b_{n,n-1}}{\lambda_n^r}\right )_{n \in \mathbf{N}}$, $\left
(\frac{b_{n,n}}{\lambda_n^r} \right )_{n \in \mathbf{N}}$ and $\left
(\frac{b_{n,n+1}}{\lambda_n^r} \right )_{n \in \mathbf{N}}$  are bounded by
$C$.
Then, for any $\psi$ in $D(|A|^r)$, $\|B\psi\| \leq \sqrt{6}C \||A|^r \psi
\|$. In particular, if $r \leq 1$ (resp. $r<1$), then  $(A,B)$ satisfies
Hypothesis \ref{HYP_1} (resp. $(A,B, \mathbf{\Phi},r)$ satisfies Hypothesis \ref{HYP_2}).
\end{prop}

\begin{pf} 
For any $\psi$ in $ D(|A|^r)$, 
\begin{eqnarray}
\lefteqn{\|B \psi\|^2=\left \|\sum_{k \in \mathbf{N}} \langle \phi_k, B\psi
\rangle \phi_k \right \|^2=\sum_{k \in \mathbf{N}} |\langle B \phi_k, \psi
\rangle |^2}\nonumber\\
&=& \sum_{k \in \mathbf{N}} |b_{k-1,k} \langle  \phi_{k-1}, \psi \rangle  + 
b_{k,k} \langle  \phi_k, \psi \rangle  +  b_{k+1,k} \langle  \phi_{k+1}, \psi
\rangle |^2\nonumber\\
&\leq&2 \sum_{k \in \mathbf{N}} |b_{k,k-1}|^2 |\langle \phi_{k-1}, \psi
\rangle|^2 + 
 |b_{k,k}|^2 |\langle \phi_{k}, \psi \rangle|^2\nonumber \\
 &&\quad \quad \quad + |b_{k,k+1}|^2 |\langle \phi_{k+1}, \psi \rangle|^2\nonumber\\
&\leq&2C^{2} \sum_{k \in \mathbf{N}} (\lambda_{k-1}^{2r} |\langle \phi_{k-1}, \psi
\rangle|^2 +
\lambda_{k}^{2r} |\langle \phi_{k}, \psi \rangle|^2\nonumber \\
 &&\quad \quad \quad + \lambda_{k+1}^{2r} |\langle \phi_{k+1}, \psi
\rangle|^2)\nonumber\\
&\leq & 6 C^2 \||A|^r \psi \|^2\label{eq:1232}
\end{eqnarray}
\end{pf}

\subsection{A toy model: the anharmonic oscillator}

Consider the system
\begin{equation}\label{EQ_harmonic_oscillator}
\mathrm{i}\frac{\partial \psi}{\partial t}(x,t)=\lbrack (-\Delta +x^2)^\alpha
+u(t) 
x^\beta\rbrack \psi(x,t),
\end{equation}
with $x$ in $\mathbf{R}$, $\psi$ in $L^2(\mathbf{R},\mathbf{C})$, $\alpha,
\beta$ in $\mathbf{N}$. 
When $\alpha=\beta=1$, (\ref{EQ_harmonic_oscillator}) is one of the most
important quantum system, it is the standard quantum harmonic oscillator submitted to a uniform electric  field.
For $\beta =1$ the system is tri-diagonal.

With our notations, $H=L^2(\mathbf{R},\mathbf{C})$, $A: \psi \mapsto -\mathrm{i}
(-\Delta +x^2)^\alpha \psi$ and $B:\psi \mapsto -\mathrm{i} x^\beta \psi$.
Operator $A$ is skew-adjoint on its domain $D(A)$, $B$ is skew-symmetric. A
Hilbert basis $\mathbf{\Phi}$ of $L^2(\mathbf{R},\mathbf{C})$ made of eigenvectors of 
$A$ is given by the sequence $(\phi_k)_{k \in \mathbf{N}}$ of the normalized Hermite
functions
$$
\phi_k:x\mapsto (-1)^k (2^k k! \sqrt{\pi})^{-1/2} e^{x^2/2}
\frac{\mathrm{d}^k}{\mathrm{d}x^k}e^{-x^2}.
$$  
For any $k$ in $\mathbf{N}$, the eigenvector $\phi_k$ is associated with the
eigenvalue $-\mathrm{i} \lambda_k=-\mathrm{i}(2k+1)^\alpha$. 

\begin{prop}\label{PRO_anharmonic_oscillator_well_posed}
If $ 2\alpha\geq \beta$ then system (\ref{EQ_harmonic_oscillator})  satifies Hypothesis 
\ref{HYP_1}.
If $2 \alpha >   \beta$ then system (\ref{EQ_harmonic_oscillator})  satifies Hypothesis 
\ref{HYP_2}.
%
%
\end{prop}
\begin{pf}
We show that the system satisfies Hypothesis~\ref{HYP_1}.\ref{ASS_B_relativ_bounded} if 
$2 \alpha\geq  \beta$ and Hypothesis~\ref{HYP_2}.\ref{ass:bdd} 
if $2 \alpha >  \beta$. The system clearly  fulfills all other hypotheses.
For any $k$ in $\mathbf{N}$,
$$x\phi_k(x)=\sqrt{\frac{k}{2}} \phi_{k-1}(x)+ \sqrt{\frac{k+1}{2}}
\phi_{k+1}(x).
$$
Iterating $\beta$ times this equality, one gets for any $k$, 
$$
|\langle x^\beta \phi_k , \psi \rangle |\leq
((k+\beta)/2)^{\beta/2}\sum_{j=-\beta}^{\beta}|\langle\phi_{k+j}, \psi \rangle| 
$$
hence, following the idea of the chain of inequalities~\eqref{eq:1232} one has
\begin{align*}
\|B\psi\|^{2} & =\sum_{k \in \mathbf{N}} |\langle B\phi_{k},\psi \rangle|^{2} \\
&\leq C \|\psi\|^{2}  + 2^{-\beta}\sum_{k > \beta} (k+\beta)^{\beta} \sum_{j=-\beta}^{\beta}|\langle\phi_{k+j}, \psi\rangle|^{2}\\
&\leq C \|\psi\|^{2}  + 2^{-\beta} \sum_{j=-\beta}^{\beta} \sum_{k > \beta} (2k+1)^{\beta} |\langle\phi_{k+j}, \psi\rangle|^{2}\\
& \leq  C \|\psi\|^{2}  + 2^{-\beta}(2\beta+1) \||A|^{\beta/(2\alpha)}\psi\|^{2},
\end{align*}
which concludes the proof.

\end{pf}
Thanks to Proposition \ref{PRO_anharmonic_oscillator_well_posed}, we can apply
Proposition \ref{PRO_continuity_norme_Ar} 
and prove the well-posedness of (\ref{EQ_harmonic_oscillator}) 
\begin{prop}
If $2\alpha = \beta$, then
(\ref{EQ_harmonic_oscillator}) is well-posed for any control $u$ with bounded variation and $L^\infty$ norm smaller than 
$\sqrt{(2\beta+1)2^{-\beta}}.$
If $2\alpha>\beta$,  then
(\ref{EQ_harmonic_oscillator}) is well-posed for any control $u$ of bounded variation.
\end{prop}

Notice that Proposition \ref{PRO_Good_Galerkin_approximtion} applies also to 
systems that are \emph{not} weakly-coupled, see \cite[Definition~1]{weakly-coupled}. For
instance, using the set 
$\{(k,k+1),k \in \mathbf{N}\}$ as a non-resonant chain of connectedness, see 
\cite[Definition~2.5]{Schrod2}, and the fact that $ (|b_{k,k+1}|^{-1})_{k\in \mathbf{N}}$ is in $\ell^1$, we get the following.
\begin{prop}
Assume that $\beta \geq 3$ odd and $\alpha>\beta/2$. Then, there exists $K=\sum_k \frac{2\pi}{|b_{k,k+1}|}>0$ such that, for any
even 
functions $\psi_0,\psi_1$  in the unit sphere of $L^2(\mathbf{R}, \mathbf{C})$, 
for any $\varepsilon>0$,
there exists a control $u_\varepsilon:[0,T_\varepsilon]\to [0,+\infty)$ such
that 
$\|\Upsilon^{u_\varepsilon}_{T_\varepsilon,0}\psi_0-\psi_1\|_{L^2}\leq 
\varepsilon$ and $\|
u_\varepsilon\|_{L^1([0,T_\varepsilon])}<K.$ 
\end{prop}
In other words, if $2\alpha>\beta\geq 3$ and $\beta$ is odd, then there is no  Good Galerkin approximation for 
(\ref{EQ_harmonic_oscillator}) in $L^2(\mathbf{R},\mathbf{C})$ in terms of
the $L^1$ 
norm of the control.
However, from Proposition \ref{PRO_Good_Galerkin_approximtion}, system
(\ref{EQ_harmonic_oscillator}) admits a sequence of Good 
Galerkin approximations  in $L^2(\mathbf{R},\mathbf{C})$ in terms of the $(L^1+TV)$ 
norm of the control.

\subsection{Rotation of a 2D molecule}\label{SEC_rotation_2D}
We consider a linear molecule whose only degree of freedom is the planar rotation, in a 
fixed plan, about its fixed center of mass.  In this model, the Schr\"{o}dinger equation 
reads
\begin{equation}\label{EQ_rotation2D}
 \mathrm{i}\frac{\partial \psi}{\partial t}=-\Delta \psi +\cos\theta \psi, \quad \theta 
 \in \Omega,
\end{equation}
$\Omega=\mathbf{R}/2\pi\mathbf{Z}$ is the unit circle endowed with the Riemannian 
structure inherited from 
$\mathbf{R}$, $H$ is the space of odd functions of $L^2(\Omega,\mathbf{C})$, 
$A=\mathrm{i}\Delta$ ($\Delta$ is the restriction to $H$ of the Laplace-Beltrami operator 
of $\Omega$) and $B:\psi \mapsto (\theta \mapsto \cos(\theta)\psi(\theta))$ is the 
multiplication by cosine. 
 
 In the Hilbert basis $\Phi=(\theta \mapsto \sin(k\theta))_{k \in \mathbf{N}}$ of $H$, $A$ 
is diagonal with diagonal $-\mathrm{i}k^2, k=1\ldots\infty$ and $B$ is tri-diagonal with 
$b_{k,k}=0, b_{k,k+1}=-\mathrm{i}/2$ for every $k$ in $\mathbf{N}$.

System \ref{EQ_rotation2D} is both tri-diagonal and weakly-coupled 
and it has been thoroughly studied (see for instance 
\cite{noiesugny-CDC} and \cite{Schrod2}). For instance, it was known that 
 (\ref{EQ_rotation2D}) admits a sequence of Good Galerkin Approximations in terms of 
$L^1$ norm of the control. More preciselyby \cite[Section IV.C]{weakly-coupled}) for every 
$\phi$ with norm $1$ in $\mathrm{span}(\phi_1,\phi_2)$, 
$$
\| X^u_{(\Phi,N)}(t,0)\phi-\pi^\Phi_N \Upsilon^u_t(\phi)\| \leq \frac{K^{N-1}}{(N-2)!}.
$$

Approximate controllability of (\ref{EQ_rotation2D}) was established in \cite{Schrod2}. 
In~\cite{periodic} is given an explicit control law to steer (\ref{EQ_rotation2D}) from $\phi_1$ 
to any neighborhood of 
$\phi_2$ using periodic functions  with  frequency $2\pi/3$. Defining $u_n:=t\mapsto 
\cos(3t)/n$ and $T^\ast=2\pi$, we have
$\displaystyle{\left | \langle \phi_2,\Upsilon^{u_n}_{nT^\ast,0}\phi_1\rangle \right |\leq 
\frac{9}{n}}$.
Since $\|B\psi\| \leq \sqrt{2}\|A\psi\|$ for every $\psi$ in $D(A)$,  Proposition 
\ref{PRO_Borne_A} implies that every control $u:[0,T]\to \{0,1\} $
with bounded variation satisfying
$ | \langle \phi_2,\Upsilon^{u}_{T,0}\phi_1\rangle  |>1-\varepsilon$  
has total variation larger than $\log(2(1-\varepsilon))/4$. 
This lower bound is rather conservative, and we will give better estimates using the 
boundedness of $B$.

For every  $u_1,u_2,t_1,t_2$ in $\mathbf{R}$, for every $\psi$ in $H$, one has
\begin{eqnarray*}
\lefteqn{\|(A+u_1B)e^{t(A+u_2B)}\psi \|}\\
&=&\|(A+u_2B)e^{t(A+u_2B)}\psi +(u_2-u_1)B  e^{t(A+u_2B)}\psi\|\\
&\leq & \|(A+u_2B)e^{t(A+u_2B)}\psi\| + \|(u_2-u_1)B  e^{t(A+u_2B)}\psi\|\\
&\leq & \|(A+u_2B) \psi \| + |u_2-u_1| \|B\| \|\psi\|
\end{eqnarray*}
For every $u_1,u_2,\ldots,u_n$ and $t_1,t_2,\ldots,t_n$ in $\mathbf{R}$, for every $\psi$ in 
the unit sphere of $H$, one shows by induction on $n$ that 
\begin{eqnarray*}
\lefteqn{\|Ae^{t_1(A+u_1B)}e^{t_2(A+u_2B)}\cdots e^{t_n(A+u_nB)}\psi \| }\\
&\leq& \|A\psi \|+\|B\|(|u_1|\!+\!|u_2-u_1|+\cdots+\!|u_n-u_{n-1}|\!+\!|u_n|)
\end{eqnarray*}
Let $k$ in $\mathbf{N}$ and $\psi$ in an $\varepsilon$-neighborhood of $\phi_k$. 
If $u$ is piecewise constant taking value in $\{0,1\}$, with $u(0)=0=\lim_{t\to\infty}u(t)$, such that $\Upsilon^u_t\phi_1=\psi$, then the number $\mathcal{N}$ of switches of  $u$ satisfies
$ \|A \psi\| \leq \|A\phi_1\| +\|B\| {\mathcal{N}}$, or
$$ 
\mathcal{N} \geq \frac{\|A\phi_{k}\| - k^{2}\varepsilon-\|A\phi_{1}\|}{\|B\|} =\frac{k^2(1-\varepsilon)-1}{\sqrt{2}}.
$$

%

\subsection{Cooling in harmonic traps}\label{SEC_cooling}

This example is inspired by \cite{jr.:1458}. The dynamics of a quantum system 
trapped in a one-dimensional parabolic potential with time varying frequency 
$\omega(t)$ is given by
\begin{equation}\label{EQ_Bose_Einstein}
\mathrm{i}\frac{\partial \psi}{\partial t}(x,t)= (-\Delta +\omega(t)
x^2 ) 
\psi(x,t),
\end{equation}
The system (\ref{EQ_Bose_Einstein}) has raised considerable attention in the
last decades (see \cite{7401448820120101} for recent developments).

Let $\lambda>0$. Defining $u(t):=\omega(t)-\lambda$, we reformulate
(\ref{EQ_Bose_Einstein}) as
\begin{equation}\label{EQ_Bose_Einstein_reform}
\mathrm{i}\frac{\partial \psi}{\partial t}(x,t)= (-\Delta + \lambda x^2 +
u(t) x^2 ) 
\psi(x,t),
\end{equation}
Note that the parity, if any, of the solutions of (\ref{EQ_Bose_Einstein}) is
preserved along the time. Hence we consider (\ref{EQ_Bose_Einstein_reform}) in
the space 
$H$ of  even functions in $L^2(\mathbf{R},\mathbf{C})$. For any $\lambda>0$, in the basis 
$\Phi=(x\mapsto \frac{1}{\lambda^{1/4}} H_{2k}(\sqrt{\lambda}x))_{k\in \mathbf{N}}$, where $H_n$ is the $n^{th}$ Hermite
functions, 
the operator $A_\lambda:=\mathrm{i}(-\Delta+\lambda^2 x^2)_{|H}$ is diagonal with diagonal 
$((2k+1)\lambda)_{k\in \mathbf{N}}$ and $B=-\mathrm{i}x^2_{|H}$ has matrix 
$[b_{j,k}]_{(j,k)\in \mathbf{N}^2}$ with $b_{j,k}=0$ if $|j-k|> 1$ and 
$b_{j,j}\sim_\infty j/\lambda$ and $b_{j,j+1}\sim_\infty j/(2\lambda)$ 
for any $j,k$ in $\mathbf{N}^2$. 

The system $(A_\lambda,B)$ is tri-diagonal and the well-posedness
of (\ref{EQ_Bose_Einstein}) follows as in Propositions
\ref{PRO_diagonal_B_A_borne}
and \ref{PRO_continuity_norme_Ar} applied to (\ref{EQ_Bose_Einstein_reform})
with $u$ any control with bounded variation and small enough.
\begin{prop}
For any even function $\psi_0$ in
$L^2(\mathbf{R},\mathbf{C})$, for any
$T>0$, for any $\alpha>0$, for any $\omega:[0,T]\to (\alpha,+\infty)$ with bounded variation, 
(\ref{EQ_Bose_Einstein}) admits a unique solution $t\mapsto
\Upsilon^\omega_{t}\psi_0$ satisfying $\Upsilon^\omega_0\psi_0=\psi_0$.
\end{prop}

\section{Conclusion}

We obtained an elementary proof of the well-posedness of bilinear 
Schr\"{o}dinger equations by adapting classical tools developed by Kato to the simple 
structure of bilinear conservative systems. 
The key ingredient of our construction is an \emph{a priori} 
upper bound on the growth of some energy functional in terms of the total variation of the 
control. 

As a consequence we prove a general method to obtain explicit  bounds on 
the number of switches of a control steering the system from a given source to a 
given target,  in the case in which the control takes value in a discrete set.
 These bounds are of importance when considering quantum
systems for which the dipolar approximation (leading to a  bilinear modeling as
in the present paper) is not valid anymore, see~\cite{Morancey} and~\cite{CDCquadratic}.

%
%
\begin{ack}
It is a pleasure for the third author to thank Ugo Boscain for inspiring discussions about 
the examples of Sections \ref{SEC_rotation_2D} and \ref{SEC_cooling}.
\end{ack}

\bibliography{biblio}             
                                                   







\end{document}